\documentclass[11pt,leqno]{amsart}
\usepackage{amssymb,verbatim,enumerate,ifthen,hyperref}
\usepackage{mathtools} 

\usepackage[mathscr]{eucal}
\usepackage[utf8]{inputenc}
\usepackage[T1]{fontenc}
\def\P{\mathbb{P}}

\def\M{\mathscr{F}}
\def\C{\mathscr{C}}

\def\M{\mathscr{M}}

\def\S{\mathscr{S}}

\long\def\comment#1{}

\newtheorem{theorem}{Theorem}[section]
\newtheorem*{theorem*}{Theorem}
\def\Thm#1#2{\ifthenelse{\equal{#1}{*}}{\begin{theorem*}#2\end{theorem*}}
             {\begin{theorem}\label{T#1}#2\end{theorem}}}
\newtheorem{Atheorem}{Theorem}

\newtheorem{proposition}[theorem]{Proposition}
\newtheorem*{proposition*}{Proposition}
\def\Prp#1#2{\ifthenelse{\equal{#1}{*}}{\begin{proposition*}#2\end{proposition*}}
{\begin{proposition}\label{P#1}#2\end{proposition}}}

\newtheorem{corollary}[theorem]{Corollary}
\newtheorem*{corollary*}{Corollary}
\def\Cor#1#2{\ifthenelse{\equal{#1}{*}}{\begin{corollary*}#2\end{corollary*}}
             {\begin{corollary}\label{C#1}#2\end{corollary}}}

\newtheorem{lemma}[theorem]{Lemma}
\newtheorem*{lemma*}{Lemma}
\def\Lem#1#2{\ifthenelse{\equal{#1}{*}}{\begin{lemma*}#2\end{lemma*}}
             {\begin{lemma}\label{L#1}#2\end{lemma}}}

\theoremstyle{definition}
\newtheorem{remark}[theorem]{Remark}
\newtheorem*{remark*}{Remark}
\def\Rem#1#2{\ifthenelse{\equal{#1}{*}}{\begin{remark}\rm #2\end{remark}}
             {\begin{remark}\label{R#1}\rm #2\end{remark}}}

\newtheorem{example}[theorem]{Example}
\newtheorem*{example*}{Example}
\def\Exa#1#2{\ifthenelse{\equal{#1}{*}}{\begin{example*}\rm #2\end{example*}}
             {\begin{example}\label{Ex#1}\rm #2\end{example}}}

\def\Eq#1#2{\ifthenelse{\equal{#1}{*}}
  {\begin{equation*}\begin{aligned}#2\end{aligned}\end{equation*}}
  {\begin{equation}\begin{aligned}\label{E#1}#2\end{aligned}\end{equation}}}

\begin{document}
\begin{flushright}
\end{flushright}
\vspace{5mm}

\date{\today}

\title[Partitioning to solve Bin Packing Problems]
{Partitioning to solve Bin Packing Problems}

\author[A. R. Goswami]{Angshuman R. Goswami}
\address[A. R. Goswami]{Department of Mathematics, University of Pannonia,
H-8200 Veszprem, Hungary}
\email{\{goswami.angshuman.robin@mik.uni-pannon.hu}


\thanks{The research of the first author was supported by the 
}

\begin{abstract}
The Bin Packing Problem involves efficiently packing items into a limited number of bins without exceeding their capacity.
In this paper, we try to answer a specific question in this field. Mathematically the combinatorial optimization problem of this classical Bin Packing can be formulated as follows\\

"Will it be possible to distribute $n$ items with sizes $a_1, \cdots ,a_n$ in total $k$ numbers of bins such that each of the bins consists of exactly $\ell$ number of items whose combined item size is $\C$ provided none of the items left after the allocation? In other words, can we achieve the optimal arrangement 
$\C k=\overset{n}{\underset{i=1}\sum} a_i$ such that $k\ell=n$ also holds under the given size restriction on bins?" \\

We propose an algorithm that determines whether a Bin Packing arrangement is possible while ensuring that each packed bin is distinct. We use tools from standard mathematical concepts such as partitioning of natural number and multiset to derive the algorithm. This algorithm is efficient, easy to implement, and ensures the desired packing in all applicable cases.\\

The following section discusses research progress, motivation, research methodologies, and various important notions and terminologies.
\end{abstract}

\maketitle

\section*{INTRODUCTION}
The Bin Packing Problem is a classic combinatorial optimization problem that arises in resource allocation in cloud computing, cutting stock problems in manufacturing, transportation logistics, and data storage optimization. The goal is to pack a given set of items into a finite number of bins in the most efficient way possible, such that the total weight or size of items in each bin does not exceed the bin's capacity. The two major variants of Bin Packing are one-dimensional Bin Packing, where items have a single attribute like weight or volume, and multi-dimensional Bin Packing, where items have multiple attributes like width, height, or depth. These can be offline, where all items are known beforehand, or online, where items arrive sequentially without prior knowledge. The Bin Packing Problem is NP-hard, meaning there is no known polynomial-time algorithm to solve it optimally for all cases. \\

The Bin Packing Problem was first studied in 1971 when Ullman introduced the concept \cite{Ullman}. An extensive study of Bin Packing Problems can be found in the Ph.D. dissertation of David S. Johnson entitled "Near-Optimal Bin Packing Algorithms" (\cite{Johnson}). Coffman\cite{Coffman} provides a comprehensive overview of early developments related to Bin Packing Problems. In the last three decades, researchers have come up with various algorithmic approaches to solve Bin Packing Problems.
The primary Heuristics algorithms such as  First Fit, Best Fit, Next Fit, and Meta Heuristics algorithms like Genetic Algorithm, Tabu Search, Ant Colony Optimization, etc. are usually used to tackle Bin Packing Problems. Moreover, several Exact and Approximation algorithms were also utilized extensively to obtain the optimal solution. Further details can be found in the papers \cite{Alvim, Borgulya, Fleszar, Fekete, Xu, Merikhipour, Scholl, Schwerin, Schoenfield} and reference therein. A recent breakthrough in this field is the integration of deep learning techniques. Additionally, reinforcement learning has been explored for Bin Packing optimization. The deep learning approach is based on maximizing space utilization and minimizing the number of bins used in order to achieve optimal packing. One can look into the papers like \cite{Wang, Zhao} to understand the recent developments properly.\\

Our objective is to propose an exact algorithm that will predict the possibility of distinct packings and provide the appropriate arrangement for it. By using the tools from classical set theory and natural number partitioning, we come up with a mathematically sound method that determines precisely whether we can obtain the packings which are different from one another or not. If no such distinct arrangements are possible, we have to do further analysis with other available algorithms to achieve the optimal result. However, in many real-world problems, distinct packing/arrangements of items is a must. For example, a mixture of drug composition focusing individual requirements, distribution of daily tasks in an office among the employees, slotting of numerous files based on data types in cloud computation, segregation of waste for proper recycling process, etc. \\

Several benchmark problems are listed in the articles \cite{Falkenauer, Delorme}. These instances can be used to compare various algorithms. In our paper, we will use the Falkanauer T-Class (\cite{Falkenauer}) to test our algorithm. Introduced by G. Falkenauer in 1996, the class serves as a standard test case to assess how effectively an algorithm can minimize the number of bins required to accommodate all items without exceeding the bin capacity. The task is to assign 60 (or 120) items in 20 (or 40) bins. The item sizes vary from 249 to 491 and the size of each of the bins is $1000$. \\

First, we broadly discuss the algorithm. We observe that this algorithm can also be used beyond Bin Packing Problems and hence we consider a general set-theoretic approach. We elaborate on the various stages involved in it and how this algorithm should be executed. Next, we describe the implementation of our algorithm in terms of Bin Packing Problems by considering the very first instance of Frankenauer T-Class, which is also referred to as Falkenauer$\_T60\_01$. The challenge is to allocate 60 items in 20 bins; each bin of size 1000. Our findings can be compared with the results in the paper \cite{Gyuri}. After this, we investigate about various pros and cons of our proposed algorithm.\\

Before presenting our algorithm, it is important to recall some concepts and terminology from classical set theory. 
A \textit{multiset} is a generalized concept of a set in mathematics and computer science, where elements are allowed to appear more than once. Unlike a standard set, which only keeps unique elements, a multiset keeps track of the multiplicity (the number of times each element appears).\\

Two multisets are considered \textit{equal} if they contain the same elements with the same multiplicities. In other words, for every element in the multisets, its count (or frequency) must be identical in both. Two multisets $\M_1$ and $\M_2$ are equal ($\M_1=\M_2$) if and only if:
\begin{enumerate}
\item Every element in $\M_1$ is in $\M_2$, and vice versa.
\item For each element $x$, the multiplicity of $x\in \M_1$ is the same as its multiplicity in $\M_2$.
\end{enumerate}

\textbf{Mathematical Notation:} 
\[
A = B \iff \forall x \left( \text{Multiplicity of } x \text{ in } A = \text{Multiplicity of } x \text{ in } B \right).
\]
The \textit{extracted set} of a multiset is a subset of it which is derived by removing the multiplicities of the elements. In simpler terms, it consists of all the unique elements present in the multiset, with each element appearing only once.\\

\section*{The Overview of the Algorithm}
Consider $n$ items of sizes $a_1,\cdots,a_n$. It is important to note that the sizes of the items may or may not be distinct. Our purpose is to arrange these $n$ items in $k$ ($1<\ell<n$) subsets such that each subset contains exactly $\ell$ number of items provided the following conditions are satisfied.

\begin{itemize}
\item The total size of $\ell$ items in each subset must be $\C.$

\item The allocation must be proper; which means after the 
specific distribution, no items are left unassigned. In other words, $k\times\ell=n$ must hold.
\end{itemize}

We propose the following algorithm, which will answer two queries. It will let us know the possibility of distinct allocation; which means any two randomly selected subsets can't have all the elements of the same sizes. If such allocation is possible, then the algorithm will generate the blueprint of that particular arrangement.\\

Here is the step-by-step description of the algorithm:
 \\
\begin{enumerate}[(1)]
\item At first we will consider the collection of all the items as the  multiset $\M=\{a_1,\cdots a_n\}$\\

\item In the next step, we compute the  extracted set $\S_{\M}=\{a_1,\cdots a_m\}$ of the multiset $\M$.\\

\item Next we will collect all the possible partitions of $\C$ that contains exactly $\ell$ elements from the set $\S_{\M}$. This newly introduced set will be denoted by the symbol $\P_{\S_{_{\M}}}^{^{\ell}}$.

$$\P_{\S_{_{\M}}}^{^{\ell}}=\{(a_1,\cdots ,a_\ell)\,\,\big| \,\,a_1+\cdots +a_\ell=\C \}.$$\\

\item After this we will consider all the subsets of $\P_{\S_{_{\M}}}^{^{\ell}}$ containing exactly $k$ elements of it. Suppose the number such subsets is $\eta.$ Obviously the value of $\eta$ depends upon the cardinality of the set $\P_{\S_{_{\M}}}^{^{\ell}}$. Mathematically we will express each one of these subsets by the symbol $\bigg(\P_{\S_{_{\M}}}^{^{\ell}}\bigg)^{k}_i$ where $i\in\{1,\cdots,\eta\}$. For instance one of such constructible subset can be as follows

$$\bigg(\P_{\S_{_{\M}}}^{^{\ell}}\bigg)^{k}_i=\{(a_1,\cdots ,a_\ell),(a_2,\cdots ,a_{\ell+1}),\cdots\cdots, (a_{k},\cdots,a_{k+\ell})\}\qquad \qquad i\in\{1,\cdots,\eta\}.$$\\

\item Following this we will construct multisets from each of the partitioning subset $\bigg(\P_{\S_{_{\M}}}^{^{\ell}}\bigg)^{k}_i$ just by spreading out all the elements from each of the tuples. In other words, each of the $\bigg(\P_{\S_{_{\M}}}^{^{\ell}}\bigg)^k_i$ corresponds to a multiset. The multiset corresponding to our previously considered $\bigg(\P_{\S_{_{\M}}}^{^{\ell}}\bigg)^{k}_{i}$ is given by
$$\M\bigg(\P_{\S_{_{\M}}}^{^{\ell}}\bigg)^{k}_{i}=\{a_1,\cdots ,a_l, a_2,\cdots ,a_{l+1},\cdots\cdots, a_{k},\cdots,a_{l+k}\}$$ \\

\item Each of these multisets  $\M\bigg(\P_{\S_{_{\M}}}^{^{\ell}}\bigg)^{k}_{i}$ $\Big(i\in\{1,\cdots,\eta\}\Big)$ is then compared with the original multiset ${\M}$ and depending upon the results we will have two conclusions.

\begin{enumerate}[(i)]
\item If $\M\bigg(\P_{\S_{_{\M}}}^{^{\ell}}\bigg)^{k}_i=\M$ holds; then the optimal allocation is successful. The partitioning set $\bigg(\P_{\S_{_{\M}}}^{^{\ell}}\bigg)^k_i$ corresponding to the multiset
$\M\bigg(\P_{\S_{_{\M}}}^{^{\ell}}\bigg)^{k}_i$ is our required solution.

\item On the other hand, if  $\M\bigg(\P_{\S_{_{\M}}}^{^{\ell}}\bigg)^{k}_i=\M$ does not hold; then there is no possibility of distinct packing.
\end{enumerate}
\end{enumerate}

\section*{Elaboration with Falkenauer T-class}
In this section, we apply our algorithm to the first benchmark problem of the Falkenauer T-class, specifically Falkenauer$\_T60\_01$. In this problem instance, we have $60$ items that vary from $251$ to $475$ in sizes whose combined size is $20,000$. In the collection, several items are of the same size also. The task is to pack these $60$ items in $20$ bins each of size $1000$ such that no items remain unpacked. Given the range of item sizes, it is evident that two items alone cannot completely fill a bin. Similarly, a bin can not accommodate four items simultaneously. Therefore, the best solution is to look for a collection of $3$ suitable items whose sizes exactly add to $1000$. Next, we attempt to fill the bins with 20 suitable triplets, ensuring that no items are left unallocated. We also observe that by this scrutinizing procedure, we transform our origin Bin Packing Problem to a $3$-partition problem.\\

So comparing this scenario with our algorithm, we have the following
\Eq{*}{
\mbox{Total number of items, $n$:}&=60\\
\mbox{Total number of bins, $k$:}&=20\\
\mbox{Total number of items in each bin, $\ell$:}&=3\\
\mbox{Total bin capacity, $\C$:}&=1000}
Now we can proceed to describe how the algorithm works in a practical sense.
\begin{itemize}
\item First consider the multiset $M$ which represents the Falkenauer$\_T60\_01$.
\Eq{*}{\M=\left\{
\begin{array}{ll}
251, 251, 252, 254, 255, 256, 257, 258, 258, 260, 260, 261, 262, 264, 265,\\ 267, 269, 270, 275, 277, 280, 282, 289, 297, 300, 302, 304, 305, 307, 308,\\ 313, 314, 319, 333, 334, 339, 340, 347, 361, 366, 369, 376, 382, 396, 396,\\ 399, 402, 403, 409, 411, 412, 423, 426, 444, 447, 462, 465, 468, 473, 475
\end{array}
\right\}.
}\\

\item In the next stage, we consider the extracted set 
\Eq{*}{
\S_{\M} =\left\{
\begin{array}{ll}
251, 252, 254, 255, 256, 257, 258, 260, 261, 262, 264, 265, 267, 269, 270, 275, 277, 280, 282,\\
289, 297, 300, 302, 304, 305, 307, 308, 313, 314, 319, 333, 334, 339, 340, 347, 361, 366,\\
369, 376, 382, 396, 399, 402, 403, 409, 411, 412, 423, 426, 444, 447, 462, 465, 468, 473, 475 
\end{array}
\right\}.
}\\

\item Now we will formulate a set that consists of all the triplets whose summation is $1000$ as follows
\Eq{*}{
\P_{\S_{_{\M}}}^{^{\ell}}=
\left\{
\begin{array}{ll}
(251, 302, 447), 
(251, 305, 444), 
(251, 340, 409), 
(251, 347, 402),\\ 
(252, 275, 473), 
(252, 280, 468), 
(252, 304, 444), 
(252, 339, 409),\\ 
(252, 366, 382), 
(254, 302, 444), 
(254, 305, 441), 
(254, 340, 406), \\ 
\vdots \vdots\vdots \vdots \vdots\vdots\vdots \vdots\vdots\vdots \vdots\vdots\vdots \vdots\vdots\vdots \vdots\vdots\vdots \vdots\vdots\vdots \vdots\vdots\vdots \vdots\vdots\vdots \vdots\vdots\vdots \vdots\vdots\vdots \vdots\vdots\vdots \vdots\vdots\vdots \vdots\vdots\vdots \vdots\vdots\vdots \vdots\vdots\vdots \vdots\vdots\vdots \vdots\vdots\vdots \vdots\vdots\vdots \vdots\vdots\vdots \vdots\vdots\vdots \vdots\vdots\vdots \vdots\vdots \vdots \vdots \vdots\vdots\vdots \vdots\vdots\vdots \vdots\vdots\vdots \vdots\vdots \vdots
\vdots\vdots\vdots \vdots\vdots\vdots \vdots\vdots \vdots
\vdots \vdots\\
\vdots \vdots\vdots \vdots \vdots\vdots\vdots \vdots\vdots\vdots \vdots\vdots\vdots \vdots\vdots\vdots \vdots\vdots\vdots \vdots\vdots\vdots \vdots\vdots\vdots \vdots\vdots\vdots \vdots\vdots\vdots \vdots\vdots\vdots \vdots\vdots\vdots \vdots\vdots\vdots \vdots\vdots\vdots \vdots\vdots\vdots \vdots\vdots\vdots \vdots\vdots\vdots \vdots\vdots\vdots \vdots\vdots\vdots \vdots\vdots\vdots \vdots\vdots\vdots \vdots\vdots\vdots \vdots\vdots \vdots \vdots \vdots\vdots\vdots \vdots\vdots\vdots \vdots\vdots\vdots \vdots\vdots \vdots
\vdots\vdots\vdots \vdots\vdots\vdots \vdots\vdots \vdots
\vdots \vdots\\
(297, 307, 396), 
(297, 334, 369), 
(300, 304, 396),
(300, 334, 366),\\ 
(300, 339, 361), 
(304, 314, 382),
(305, 313, 382), 
(305, 319, 376),\\ 
(305, 334, 361),
(313, 340, 347), 
(314, 339, 347), 
(319, 334, 347)
\end{array}
\right\}.
}
The set $\P_{\S_{_{\M}}}^{^{\ell}}$ consists of $99$ triplets.
\\

\item Now we will consider all subsets of $\P_{\S_{_{\M}}}^{^{\ell}}$ of size $20$. Total number of such subsets will be $\binom{99}{20}=428,786,696,323,047,746,376.$ These subsets can be listed as below
\Eq{*}{
\bigg(\P_{\S_{_{\M}}}^{^{\ell}}\bigg)^{k}_1=
\left\{
\begin{array}{ll}
(251, 302, 447), 
(251, 305, 444), 
(251, 340, 409), 
(251, 347, 402), 
(252, 275, 473),\\ 
(252, 280, 468), 
(252, 304, 444), 
(252, 339, 409), 
(252, 366, 382), 
(254, 302, 444),\\
(254,334,412),
(254,347,399),
(255,270,475),
(255,277,468),
(255,280,465),\\
(255,319,426),
(255,333,412),
(255,334,411),
(255,369,376),
(256,269,475),
\end{array}
\right\};
}
\Eq{*}{
\bigg(\P_{\S_{_{\M}}}^{^{\ell}}\bigg)^{k}_2=
\left\{
\begin{array}{ll}
(251, 302, 447), 
(251, 305, 444), 
(251, 340, 409), 
(251, 347, 402), 
(252, 275, 473),\\ 
(252, 280, 468), 
(252, 304, 444), 
(252, 339, 409), 
(252, 366, 382), 
(254, 305, 441),\\
(254,334,412),
(254,347,399),
(255,270,475),
(255,277,468),
(255,280,465),\\
(255,319,426),
(255,333,412),
(255,334,411),
(255,369,376),
(256,282,462)
\end{array}
\right\};
}
\Eq{*}{
\vdots \vdots\vdots \vdots \vdots\vdots\vdots \vdots\vdots\vdots \vdots\vdots\vdots \vdots\vdots\vdots \vdots\vdots\vdots \vdots\vdots\vdots \vdots\vdots\vdots \vdots\vdots\vdots \vdots\vdots\vdots \vdots\vdots\vdots \vdots\vdots\vdots \vdots\vdots\vdots \vdots\vdots\vdots \vdots\vdots\vdots \vdots\vdots\vdots \vdots\vdots\vdots \vdots\vdots\vdots \vdots\vdots\vdots \vdots\vdots\vdots \vdots\vdots\vdots \vdots\vdots\vdots \vdots\vdots \vdots \vdots \vdots\vdots\vdots \vdots\vdots\vdots \vdots\vdots\vdots \vdots\vdots \vdots
\vdots\vdots\vdots \vdots\vdots\vdots \vdots\vdots \vdots
\vdots\vdots\vdots \vdots\vdots\vdots \vdots\vdots\vdots \vdots\vdots\vdots \vdots\vdots\vdots \vdots\vdots\vdots \vdots\vdots\vdots \vdots\vdots\vdots \vdots\vdots\vdots \vdots\vdots\vdots \vdots\vdots \vdots \vdots \vdots\vdots\vdots \vdots\vdots\vdots \vdots\vdots\vdots \vdots\vdots \vdots
\vdots\vdots\vdots \vdots\vdots\vdots \vdots\vdots \vdots\\
\vdots \vdots\vdots \vdots \vdots\vdots\vdots \vdots\vdots\vdots \vdots\vdots\vdots \vdots\vdots\vdots \vdots\vdots\vdots \vdots\vdots\vdots \vdots\vdots\vdots \vdots\vdots\vdots \vdots\vdots\vdots \vdots\vdots\vdots \vdots\vdots\vdots \vdots\vdots\vdots \vdots\vdots\vdots \vdots\vdots\vdots \vdots\vdots\vdots \vdots\vdots\vdots \vdots\vdots\vdots \vdots\vdots\vdots \vdots\vdots\vdots \vdots\vdots\vdots \vdots\vdots\vdots \vdots\vdots \vdots \vdots \vdots\vdots\vdots \vdots\vdots\vdots \vdots\vdots\vdots \vdots\vdots \vdots
\vdots\vdots\vdots \vdots\vdots\vdots \vdots\vdots \vdots
\vdots\vdots\vdots \vdots\vdots\vdots \vdots\vdots\vdots \vdots\vdots\vdots \vdots\vdots\vdots \vdots\vdots\vdots \vdots\vdots\vdots \vdots\vdots\vdots \vdots\vdots\vdots \vdots\vdots\vdots \vdots\vdots \vdots \vdots \vdots\vdots\vdots \vdots\vdots\vdots \vdots\vdots\vdots \vdots\vdots \vdots
\vdots\vdots\vdots \vdots\vdots\vdots \vdots\vdots \vdots
}
\Eq{*}{
\bigg(\P_{\S_{_{\M}}}^{^{\ell}}\bigg)^{k}_{\binom{99}{20}}=
\left\{
\begin{array}{ll} 
(280,308,412),
(282,307,411),
(282,319,399),
(289,300,411),
(289,302,409),\\
(289,308,403),
(297,300,403),
(297,304,399),
(297,307,396),
(297,334,369),\\
(300, 304, 396),
(300, 334, 366),
(300, 339, 361), 
(304, 314, 382),
(305, 313, 382),\\ 
(305, 319, 376), 
(305, 334, 361),
(313, 340, 347), 
(314, 339, 347), 
(319, 334, 347)
\end{array}
\right\}.
}
\\
\item Based on these subsets $\bigg(\P_{\S_{_{\M}}}^{^{\ell}}\bigg)^{k}_i$, we form the corresponding multisets $\M\bigg(\P_{\S_{_{\M}}}^{^{\ell}}\bigg)^{k}_i$ for all \newline 
$i\in\Big\{1,\cdots, \binom{99}{20} \Big\}$ as described in the fifth point of the previous section. Each of these multisets contains various sizes of $3\times 20=60$ items.\\
For instance the multiset corresponding to $\bigg(\P_{\S_{_{\M}}}^{^{\ell}}\bigg)^{k}_1$ is
\Eq{*}{
\M\bigg(\P_{\S_{_{\M}}}^{^{\ell}}\bigg)^{k}_1=
\left\{
\begin{array}{ll}
251, 302, 447, 
251, 305, 444, 
251, 340, 409, 
251, 347, 402, 
252, 275, 473,\\ 
252, 280, 468, 
252, 304, 444, 
252, 339, 409, 
252, 366, 382, 
254, 302, 444,\\
254,334,412,
254,347,399,
255,270,475,
255,277,468,
255,280,465,\\
255,319,426,
255,333,412,
255,334,411,
255,369,376,
256,269,475,
\end{array}
\right\}.
}

Similarly the multiset corresponding to $\bigg(\P_{\S_{_{\M}}}^{^{\ell}}\bigg)^{k}_2$ is
\Eq{*}{
\M\bigg(\P_{\S_{_{\M}}}^{^{\ell}}\bigg)^{k}_2=
\left\{
\begin{array}{ll}
251, 302, 447, 
251, 305, 444, 
251, 340, 409, 
251, 347, 402, 
252, 275, 473,\\ 
252, 280, 468, 
252, 304, 444, 
252, 339, 409, 
252, 366, 382, 
254, 305, 441,\\
254,334,412,
254,347,399,
255,270,475,
255,277,468,
255,280,465,\\
255,319,426,
255,333,412,
255,334,411,
255,369,376,
256,282,462
\end{array}
\right\}.}

And so on.\\
\item In the final step we compare each one of these multisets 
$\M\bigg(\P_{\S_{_{\M}}}^{^{\ell}}\bigg)^{k}_i$ for all 
$i\in\Big\{1,\cdots,\binom{99}{20} \Big\}$ with our original multiset $\M$; until we find the equality. In the case of Falkenauer$\_T60\_01$, it turns out that one of the multisets $\M\bigg(\P_{\S_{_{\M}}}^{^{k}}\bigg)^{\ell}_j$ $\Big(1\leq j\leq 15,579,278,510,796\Big)$ satisfies this equality requirement. The corresponding set of triplets is as follows
\Eq{*}{
\bigg(\P_{\S_{_{\M}}}^{^{k}}\bigg)^{\ell}_j
\left\{
\begin{array}{ll}
(251,302,447), (251,305,444), (252,339,409), (254,347,399),(255,280,465),\\
(256,269,475), (257,361,382), (258,319,423), (258,366,376),(260,267,473),\\
(260,314,426),(261,277,462),(262,270,468), (264, 340,396), (265,333,402),\\ 
(275,313,412), (282,307,411), (289,308,403), (297,334,369)
(300,304,396)
\end{array}
\right\}.}
\end{itemize}
This provides the optimal solution for the first instance of Falkenauer T class.
\\

The algorithm is grounded in strong mathematical theory, using concepts such as set partitioning and multisets. This solid foundation allows it to function as an exact algorithm, ensuring that any solution it produces strictly meets the criteria of the bin packing problem. It checks precisely whether each bin has a different combination of items, which is useful in applications like file organization or task distribution. It performs well on small problems and provides accurate, reliable results.\\

However, the algorithm has some computational limitations. It requires both the number of items per bin and the total bin size to be fixed. Another issue is that it can become time-consuming for large input sizes, as it attempts all possible groupings. This slows down the process and limits quick solutions. Nonetheless, with ongoing advances in computational power, such limitations may gradually be overcome.

\begin{thebibliography}{1}
\bibitem{Ullman}
Ullman, Jeffrey D. 
\newblock{"The Performance of a Memory Allocation Algorithm."} Technical Report 100, Princeton University, (1971).

\bibitem{Johnson}
Johnson, David S. 
\newblock{"Near-optimal bin packing algorithms."}
Diss. Massachusetts Institute of Technology, 1973.

\bibitem{Coffman}
Coffman Jr, Edward G., Michael R. Garey, and David S. Johnson. \newblock{"Approximation algorithms for bin-packing—an updated survey."}
\newblock{Algorithm design for computer system design.} Vienna: Springer Vienna, 1984. 49-106.

\bibitem{Edward}
Coffman Jr, Edward G., Csirik, J., Galambos, G., Martello, S., and Vigo, D. 
\newblock{"Bin Packing Approximation Algorithms: Survey and Classification."}
\newblock{Handbook of combinatorial optimization.}
 Springer, 2013. 455-531.
 


\bibitem{Alvim}
Alvim, A. C., Ribeiro, C. C., Glover, F., and Aloise, D. J. 
\newblock{"A hybrid improvement heuristic for the one-dimensional bin packing problem."}
\newblock{Journal of heuristics} 10 (2004): 205-229.

\bibitem{Borgulya}
Borgulya, Istvan. 
\newblock{"A hybrid evolutionary algorithm for the offline Bin Packing Problem."} 
\newblock{Central European Journal of Operations Research} 29.2 (2021): 425-445.

\bibitem{Fleszar}
Fleszar, Krzysztof, and Christoforos Charalambous. 
\newblock{"Average-weight-controlled bin-oriented heuristics for the one-dimensional bin-packing problem."} 
\newblock{European Journal of Operational Research} 210.2 (2011): 176-184.

\bibitem{Fekete}
Fekete, Sándor P., and Jörg Schepers.
\newblock{"New classes of fast lower bounds for bin packing problems."}
\newblock{Mathematical programming} 91 (2001): 11-31.

\bibitem{Xu}
Hongjie Xu, Yunzhuang Shen, Yuan Sun, Xiaodong Li. 
\newblock{"Machine Learning-Enhanced Ant Colony Optimization for Column Generation."} 
\newblock{Proceedings of the Genetic and Evolutionary Computation Conference.} 2024.

\bibitem{Merikhipour}
Merikhipour, Mahsa, Shayan Khanmohammadidoustani, and Mohammadamin Abbasi. 
\newblock{"Transportation mode detection through spatial attention-based transductive long short-term memory and off-policy feature selection."}
\newblock{Expert Systems with Applications} 267 (2025): 126196.

\bibitem{Scholl}
Scholl, Armin, Robert Klein, and Christian Jürgens. 
\newblock{"Bison: A fast hybrid procedure for exactly solving the one-dimensional bin packing problem."}
\newblock{Computers \& Operations Research} 24.7 (1997): 627-645.

\bibitem{Schwerin}
Schwerin, Petra, and Gerhard Wäscher. 
\newblock{"The bin-packing problem: A problem generator and some numerical experiments with FFD packing and MTP."}
\newblock{International transactions in operational research} 4.5-6 (1997): 377-389.

\bibitem{Schoenfield}
Schoenfield, Jon E. 
\newblock{"Fast, exact solution of open bin packing problems without linear programming."}
\newblock{Draft}, US Army Space and Missile Defense Command, Huntsville, Alabama, USA (2002).

\bibitem{Wang}
Wang, B., Lin, Z., Kong, W., and Dong, H. 
\newblock{"Bin Packing Optimization via Deep Reinforcement Learning."}
\newblock{ IEEE Robotics and Automation Letters} (2025).

\bibitem{Zhao}
Zhao, H., Zhu, C., Xu, X., Huang, H., and Xu, K. 
\newblock{"Learning practically feasible policies for online 3D bin packing."}
\newblock{Science China Information Sciences}
 65.1 (2022): 112105.

\bibitem{Falkenauer}
Falkenauer, Emanuel. 
\newblock{"A hybrid grouping genetic algorithm for bin packing."}
\newblock{Journal of Heuristics} 2 (1996): 5-30.


\bibitem{Delorme}
Delorme, Maxence, Manuel Iori, and Silvano Martello. 
\newblock{"BPPLIB: a library for bin packing and cutting stock problems."}
\newblock{Optimization Letters} 12 (2018): 235-250.

\bibitem{Gyuri}
Dósa, György, András Éles, Angshuman Robin Goswami, István Szalkai, and Zsolt Tuza. \newblock{"Solution of Bin Packing Instances in Falkenauer T Class: Not So Hard."}
\newblock{Algorithms 18, no. 2 (2025): 115.}
\end{thebibliography}

\end{document}